\documentclass[11pt]{article}
\usepackage{amssymb}
\usepackage{color}

\newcommand{\bN}{\mathbf{N}}

\newcommand{\bK}{\mathbf{K}}
\newcommand{\ord}{\mbox{\rm ord }}

\newlength{\szer}

\newtheorem{defi}{Definition}[section]
\newtheorem{nota}[defi]{Remark}

\newtheorem{teorema}[defi]{Theorem}
\newtheorem{prop}[defi]{Proposition}
\newtheorem{lema}[defi]{Lemma}
\newtheorem{coro}[defi]{Corollary}

\newenvironment{proof}[1][Proof]{\textbf{#1.} }{\
\rule{0.5em}{0.5em}}

\begin{document}
\title{On the intersection multiplicity of plane branches
\footnotetext{
     \noindent   \begin{minipage}[t]{4in}
       {\small
       2000 {\it Mathematics Subject Classification:\/} Primary 32S05;
       Secondary 14H20.\\
       Key words and phrases: intersection multiplicity, semigroup associated with a branch, key polynomials, logarithmic distance, van der Kulk theorem.\\
       The first-named author was partially supported by the Spanish Project
    MTM 2016-80659-P.}
       \end{minipage}}}

\author{Evelia R.\ Garc\'{\i}a Barroso and Arkadiusz P\l oski}

\maketitle

\begin{abstract}
\noindent We prove an intersection formula for two plane branches in terms of their semigroups and key polynomials. Then we provide a strong version of Bayer's theorem on the set of intersection numbers of two branches and apply it to the logarithmic distance in the space of branches.
\end{abstract}

\section{Introduction}
\label{intro}
\noindent Let $f, g$ be irreducible power series in $\bK[[x,y]]$, where $\bK$ is an algebraically closed field. The intersection multiplicity $i_0(f,g)$ of branches $\{f=0\}$ and $\{g=0\}$ is a basic notion of the local geometry of plane algebraic curves. The classical formula for $i_0(f,g)$ allows to calculate the intersection multiplicity in terms of Puiseux parametrisations (see \cite{van der Waerden}, \cite{Hefez}, when char $\bK=0$). If $\bK$ has a positive characteristic  a similar result can be obtained by using the Hamburger-Noether expansions (see \cite{Ancochea}, \cite{Campillo}, \cite{Russell}). The aim of this note is to prove a formula for $i_0(f,g)$ in terms of semigroups $\Gamma(f)$ and $\Gamma(g)$ associated with $f$ and $g$ and key polynomials $f_i$ and $g_j$ which define (in generic coordinates) the maximum contact curves $\{f_i=0\}$ and $\{g_j=0\}$. We impose no condition on the characteristic of $\bK$. We will use the basic notions and theorems of the theory of plane branches developed in our article \cite{GB-P} without resorting to the Hamburger-Noether expansions.

\medskip

\noindent  In Section \ref{preliminares} we recall the main properties of the semigroup and key polynomials associated with a branch. In Sections \ref{main theorem} and \ref{proof} we present the main result  (Theorem \ref{igualdad polinomios}) and its proof. Then we give in Section \ref{Kulk} an application of the main result to polynomial automorphisms of the affine plane (Theorem \ref{van der Kulk}). In Section \ref{sBayer} we prove a strong version  (due to Hefez in characteristic $0$) of Bayer's theorem on the set of intersection numbers of two branches (Theorem \ref{Bayer}). Section \ref{log dist} is devoted to a short proof of a property of the logarithmic distance $d(f,g)=\frac{i_0(f,g)}{\ord f \ord g}$ (Theorem \ref{Carr})  discovered by Ab\'{\i}o et al. in \cite{Alberich} in the case of characteristic $0$.

\section{Preliminaries}
\label{preliminares}

\noindent In this note we use the basic notions and theorems of the theory of plane branches explained in \cite{GB-P}.

\medskip

\noindent Let $\bK$ be an algebraically closed field of arbitrary characteristic. For any power series $f,g\in \bK[[x,y]]$ we define the
 {\em intersection multiplicity or intersection number} $i_0(f,g)$ by putting
\[i_0(f,g)=\mathrm{dim}_{\bK}\bK[[x,y]]/(f,g), \]

\noindent where $(f,g)$ is the ideal of $\bK[[x,y]]$ generated by $f$ and $g$. If $f,g$ are non-zero power series without constant term then $i_0(f,g)<+\infty$ if and only if 
$f$ and $g$ are coprime.

\medskip

\noindent Let $f\in \bK[[x,y]]$ be an irreducible power series. By definition, the branch $\{f=0\}$ is the ideal generated by $f$ in $\bK[[x,y]]$.

\noindent For any branch $\{f=0\}$ we put

\[\Gamma(f)=\{i_0(f,g)\;:\; g\; \hbox{\rm runs over all power series such
that } g\not\equiv 0 \;\hbox{\rm (mod $f$)}  \}.\]

\noindent Then $\Gamma(f)$ is a semigroup. We call $\Gamma(f)$ the {\em semigroup associated with the branch} $\{f=0\}$. 
\medskip

\noindent Two branches $\{f=0\}$ and $\{g=0\}$ are {\em equisingular} if and only if $\Gamma(f)=\Gamma(g)$. The branch $\{f=0\}$ is {\em non-singular} (that is of multiplicity $1$) if and only if $\Gamma(f)=\bN$. We have $\min (\Gamma(f)\backslash\{0\})=\ord f$.

\medskip

\noindent Let $n>0$ be an integer. A sequence of positive integers $(v_0,\ldots, v_h)$ is said to be an $n$-{\em characteristic sequence} if $v_0=n$ and if the following two conditions are fulfilled:

\begin{enumerate}
\item [{\rm (char 1)}] Let $e_k=\gcd(v_0,\ldots,v_k)$ for $0\leq k\leq h$. Then $n=e_0>e_1>\cdots>e_h=1$. 
\item [{\rm (char 2)}]  $e_{k-1}v_k< e_kv_{k+1}$ for $1\leq k\leq h-1$.
\end{enumerate}

\noindent Let $n_k=\frac{e_{k-1}}{e_k}$ for $1\leq k\leq h$.

\medskip

\noindent Conditions (char 1) and (char 2) imply  {\em B\'ezout's relation}:
\[
n_kv_k=a_0v_0+a_1v_1+\cdots+a_{k-1}v_{k-1},
\]

\noindent where $a_0>0$, $0\leq a_i<n_i$ for $1\leq i\leq k$ are integers.

\medskip

\noindent  Let $f=f(x,y)$ be an irreducible power series. Suppose that  $n=i_0(f,x)=\ord f(0,y)< +\infty$. Then  the $n$-{\em minimal system of generators} of $\Gamma(f)$, defined by the conditions:

\begin{itemize}
\item[\rm{(gen 1)}] $v_0=n$, $v_k$ is the smallest element of $\Gamma(f)$ which does not belong to $v_0 \bN+\cdots +v_{k-1} \bN$.
\item[\rm{(gen 2)}] $v_0 \bN+\cdots +v_{h} \bN=\Gamma(f)$
\end{itemize}

\noindent is an $n$-characteristic sequence. We will call it the characteristic of $\{f=0\}$ and write
 $\overline{\hbox{\rm char}}_x f=(v_0, \ldots, v_h)$. If $(v_0, \ldots, v_h)$ is the $n$-sequence of generators of $\Gamma(f)$ then the number $c=\sum_{k=1}^h (n_k-1)v_k-v_0+1$ is the {\em conductor} of $\Gamma(f)$, that is, $c+N\in \Gamma(f)$ for $N\geq 0$ and $c-1\not\in \Gamma(f)$.

\medskip

\noindent There exists a sequence of monic polynomials $f_0, f_1, \ldots, f_{h-1}\in \bK[[x]][y]$ such that $\deg_y f_k=\frac{n}{e_k}$ and $i_0(f,f_k)=v_{k+1}$ for $k=0, \ldots, h-1$.

\medskip

\noindent Let $f_h\in \bK[[x]][y]$ be the distinguished polynomial associated with $f$. Then $\deg_y f_h=\frac{n}{e_h}=n$ and $i_0(f,f_h)=+\infty$. We put $v_{h+1}=
+\infty$.

\medskip

\noindent The polynomials $f_0, f_1, \ldots, f_{h}\in \bK[[x]][y]$ are called {\em key polynomials} of $f$. They are not uniquely determined by $f$. Recall that for any $n$-characteristic sequence $(v_0,v_1,\ldots,v_h)$ there exists an irreducible series $f$ such that $\Gamma(f)=v_0\bN+\cdots v_h\bN$ and $i_0(f,x)=v_0$ (see \cite[Theorem 6.5]{GB-P}).

\medskip

\noindent The basic properties of key polynomials are:

\begin{enumerate}
\item [\rm{(key 1)}] A key polynomial  $f_k$ is a distinguished, irreducible polynomial and $\overline{\hbox{\rm char}}_x f_k=\left(\frac{v_0}{e_k}, \ldots, \frac{v_k}{e_k}\right)$,  consequently $\deg_y f_k=\frac{v_0}{e_k}$.

\item [\rm{(key 2)}] Let $(v_0, \ldots, v_k)$ be an $n$-characteristic sequence. Let  $f_k$ be an irreducible distinguished polynomial  such that $\overline{\hbox{\rm char}}_x f_k=\left(\frac{v_0}{e_k}, \ldots, \frac{v_k}{e_k}\right)$. Let $f_0,f_1,\ldots,f_{k-1}$ be a sequence of key polynomials of $f_k$. 

\noindent Put $f_i=f_{i-1}^{n_i}+\xi_ix^{a_{i,0}}f_0^{a_{i,1}}\cdots f_{i-2}^{a_{i,i-1}}$ for $i\in\{k+1,\ldots, h\}$, where $n_iv_i=a_{i,0}v_0+\cdots+a_{i,i-1}v_{i-1}$ is a B\'ezout's relation and $\xi_i\in \bK\backslash\{0\}$. Then $\overline{\hbox{\rm char}}_x f_h=\left(v_0, \ldots, v_h\right)$ and $f_0,\ldots, f_{h-1}$ are key polynomials  of $f_h$.

\item  [\rm{(key 3)}]  If $g\in \bK[[x]][y]$ is a monic polynomial such that $\deg_y g=\frac{n}{e_k}$ then $i_0(f,g)\leq i_0(f,f_k)=v_{k+1}$.

\end{enumerate}

\noindent For the proofs of (key 1) and (key 2) we refer the reader to \cite[Proposition 4.2 and Theorem 6.1]{GB-P}. For the proof of (key 3) see  \cite[Lemma 3.12]{GB-P}.

\medskip

\noindent The following result is a local version of the Abhyankar-Moh result:

\medskip

\noindent {\bf Abhyankar-Moh irreducibility criterion.-}
\label{AM-irr}
Let $f(x,y)\in \bK[[x,y]]$ be an irreducible power series such that $n=i_0(f,x)<+\infty$ and let $\Gamma(f)=\bN v_0+\cdots \bN v_h$, where $v_0=n$. If $g(x,y)\in \bK[[x,y]]$ is a power series such that $i_0(g,x)=n$ and  $i_0(f,g)>e_{h-1}v_h$ then $g$ is irreducible and $\Gamma(g)=\Gamma(f)$.

\medskip

\noindent The proof of the above criterion is given in \cite[Corollary 5.8]{GB-P}.

\section{Main result}

\label{main theorem}

\noindent Let $\{f=0\}$ and $\{g=0\}$ be two branches different from $\{x=0\}$. Let $\overline{\hbox{\rm char}}_x f=(v_0, \ldots, v_h)$, where $v_0=n=i_0(f,x)$ and $\overline{\hbox{\rm char}}_x g=(v'_0, \ldots, v'_{h'})$, where $v'_0=n'=i_0(g,x)$. We denote by $f_0, f_1, \ldots, f_{h}$ and $g_0, g_1, \ldots, g_{h'}$ a sequence of key polynomials of $f$ and $g$, respectively.

\medskip

\begin{teorema}(Intersection formula)
\label{igualdad polinomios} 
With the assumptions and notations introduced above there is an integer $0<k\leq \min\{h,h'\}+1$ such that

\begin{enumerate}
\item[\rm{(a)}] $\frac{v_i}{n}= \frac{v'_i}{n'}$ for all $i<k$.
\item[\rm{(b)}]$i_0(f,g)\leq \inf \{e'_{k-1}v_k,e_{k-1}v'_k\}$.
\item[\rm{(c)}] If $i_0(f,g)< \inf \{e'_{k-1}v_k,e_{k-1}v'_k\}$ then $i_0(f,g)=e_{k-1}e'_{k-1}i_0(f_{k-1},g_{k-1})$.
\item[\rm{(d)}] Suppose that $k>1$. Then $i_0(f,g)> \inf \{e'_{k-2}v_{k-1},e_{k-2}v'_{k-1}\}$.
\end{enumerate}
\noindent Moreover $f_0, f_1, \ldots, f_{k-2}$ are the first $k-1$  polynomials of a sequence of key polynomials of $g$ and $g_0, g_1, \ldots, g_{k-2}$ are the first $k-1$ polynomials of a sequence of key polynomials of $f$.
\end{teorema}

\begin{nota}
\label{nota}
From the first part of Theorem \ref{igualdad polinomios}  it follows that 
$\frac{n}{e_i}=\frac{n'}{e'_i}$ for $i<k$. In fact, 

\begin{eqnarray*}
ne'_i&=&n\gcd(v'_0, \ldots, v'_i)=
\gcd(nv'_0, \ldots, nv'_i)=\gcd(n'v_0, \ldots, n'v_i) \\
& = &n'\gcd(v_0, \ldots, v_i)=n'e_i.
\end{eqnarray*}

\noindent Consequently $n_i=n'_i$ and
$e'_{i-1}v_{i}=e_{i-1}v'_{i}$ for $0<i<k$.
\end{nota}

\noindent The proof of Theorem \ref{igualdad polinomios} is given in Section \ref{proof}. Observe that the integer $k>0$ is the smallest integer such that condition (b) of Theorem  \ref{igualdad polinomios} holds.

\begin{coro} (see \cite[Lemma 1.7]{Delgado})
\label{c:1}
Let $k>0$ be the minimum integer such that 
\[i_0(f,g)\leq \inf \{e'_{k-1}v_k,e_{k-1}v'_k\}.\]

\noindent Then 

\begin{enumerate}
\item [\rm{(1)}] $\frac{v_i}{n}= \frac{v'_i}{n'}$ for all $i<k$.
\item [\rm{(2)}] If $i_0(f,g)< \inf \{e'_{k-1}v_k,e_{k-1}v'_k\}$ then $i_0(f,g)\equiv 0$ {\rm (mod }$e_{k-1}e'_{k-1}${\rm )}.
\end{enumerate}
\end{coro}

\noindent F. Delgado in \cite[Section 3]{Delgado1} and \cite[pp. 335-336]{Delgado} computed the integer $k>0$ in terms of the Hamburger-Noether expansions of $f$ and $g$. In the sequel we do not need any additional information about the number $k$.
\section{Proof of Theorem \ref{igualdad polinomios} }
\label{proof}

\noindent For any branches $\{f=0\}$ and $\{g=0\}$ different from the branch $\{x=0\}$ we put

\[
d_x(f,g)=\frac{i_0(f,g)}{i_0(f,x)i_0(g,x)}.
\]

\noindent The function $d_x$ satisfies the {\em Strong Triangle Inequality (STI)}: for any branches $\{f=0\}$, $\{g=0\}$ and $\{h=0\}$ different from $\{x=0\}$:
\[
d_x(f,g)\geq \inf \{d_x(f,h),d_x(g,h)\},
\]
which is equivalent to the following: {\em at least two of the numbers $d_x(f,g),d_x(f,h),$ $d_x(g,h)$ are equal and the third is not smaller than the other two} (see \cite[Section 2, Theorem 2.8]{GB-P}).

\medskip

\begin{lema}
\label{lisa}
If $n=1$ then $i_0(f,g)=\inf \{e_0'i_0(f_0,g_0), v_1'\}$.
\end{lema}

\noindent \begin{proof}
If $n=1$ then $e_0=1$, $h=0$ and the only possible value for $k$, $0<k\leq \min(h+1,h'+1)$ is $k=1$. Note that $d_x(f,g)=\frac{i_0(f,g)}{v'_0}$, $d_x(f,g_0)=i_0(f,g_0)$, $d_x(g,g_0)=\frac{v'_1}{v'_0}\not\in \bN$ (if $n'=1$ then $v'_1=+\infty$ and $d_x(g,g_0)=+\infty$). Therefore $d_x(f,g_0)\neq d_x(g,g_0)$ and $d_x(f,g)=\inf \{d_x(f,g_0),d_x(g,g_0)\}$, which is equivalent to $i_0(f,g)=\inf\{e_0'i_0(f,g_0), v_1'\}=\inf \{e_0'i_0(f_0,g_0), v_1'\}.$
\end{proof}

\medskip

\noindent Lemma \ref{lisa} implies Theorem \ref{igualdad polinomios} in the case $n=1$. In the sequel we assume $n>1$.

\medskip

\noindent Let $f_0,\ldots,f_h$ be a sequence of key polynomials of $f$. Then $d_x(f,f_{k-1})=\frac{e_{k-1}v_k}{n^2}$. Thus the sequence $d_x(f,f_{k-1})$, for $k=1,\ldots,h$ is strictly increasing.

\begin{prop}
\label{ppp}
Suppose that $d_x(f,g)>d_x(f,f_{k-1})$ for an integer $k\in \{1, \ldots,h\}$. Then $k\leq h'$, $d_x(f,f_{i-1})=d_x(g,g_{i-1})$ for $i=1,\ldots,k$ and $f_0,\ldots,f_{k-1}$  are the first $k-1$ polynomials of a sequence of  key polynomials of $g$.
\end{prop}
\noindent \begin{proof} See \cite[Theorem 5.2]{GB-P}. \end{proof}

\begin{lema}
\label{qqq}
Let $k\leq  \min\{h,h'\}+1$. Suppose that $d_x(f,f_{i-1})=d_x(g,g_{i-1})$ for $0< i<k$ and $d_x(f,f_{k-1})\neq d_x(g,g_{k-1})$. Then \[d_x(f,g)\leq \inf \{d_x(f,f_{k-1}),d_x(g,g_{k-1})\}.\]
\end{lema}
\noindent \begin{proof}
From $d_x(f,f_{i-1})=d_x(g,g_{i-1})$ for $0< i<k$ we get $\frac{v_i}{n}= \frac{v'_i}{n'}$  and
$\frac{e_i}{n}= \frac{e'_i}{n'}$ for $0<i<k$. Thus $\deg_y f_{k-1}=\frac{n}{e_{k-1}}=\frac{n'}{e'_{k-1}}=\deg_y g_{k-1}$. We may assume that 
$d_x(g,g_{k-1})< d_x(f,f_{k-1})$. Since  $\deg_y f_{k-1}=\deg_y g_{k-1}$ we get, applying 
(key 3)  to $g$ that $i_0(g,f_{k-1})\leq i_0(g,g_{k-1})$ and consequently  $d_x(g,f_{k-1})\leq d_x(g,g_{k-1})$.
Thus $d_x(g,f_{k-1})< d_x(f,f_{k-1})$ and by the STI $d_x(g,f_{k-1})=d_x(f,g)$. Therefore $d_x(f,g)=d_x(g,f_{k-1})\leq d_x(g,g_{k-1})=\inf\{d_x(f,f_{k-1}),d_x(g,g_{k-1})\}.$
\end{proof}

\begin{prop}
\label{rrr}
Let $0<k\leq h+1$ be the smallest integer such that $d_x(f,g)\leq d_x(f,f_{k-1})$. Then $k\leq h'+1$ and $d_x(f,g)\leq d_x(g,g_{k-1})$.
\end{prop}

\noindent \begin{proof} First we suppose that $k=1$. Then the inequality $k\leq h'+1$ is obvious and the proposition follows from Lemma \ref{qqq}.

\medskip

\noindent Suppose now that $k>1$. By definition of $k$ we have $d_x(f,g)>d_x(f,f_{k-2})$. Then by Proposition  \ref{ppp} $k\leq h'+1$ and $d_x(f,f_{i-1})=d_x(g,g_{i-1})$ for $i=1,\ldots,k-1$.  
If $d_x(f,f_{k-1})=d_x(g,g_{k-1})$ then the proposition is obvious. If $d_x(f,f_{k-1})\neq d_x(g,g_{k-1})$ then we use Lemma \ref{qqq}.
\end{proof}

\medskip

\noindent{\bf Proof of Theorem \ref{igualdad polinomios} } Recall that $n>1$. The assertions of the theorem can be rewritten in the following form:

\begin{enumerate}
\item[\rm{(a')}] if $k>1$ then $d_x(f,f_{i-1})=d_x(g,g_{i-1})$ for all $0<i<k$.
\item[\rm{(b')}]$d_x(f,g)\leq \inf \{d_x(f,f_{k-1}), d_x(g,g_{k-1})\}$.
\item[\rm{(c')}] If $d_x(f,g)<  \inf \{d_x(f,f_{k-1}), d_x(g,g_{k-1})\}$ then $d_x(f,g)=d_x(f_{k-1},g_{k-1})$.
\item[\rm{(d')}] If $k>1$ then $d_x(f,g)>  \inf \{d_x(f,f_{k-2}), d_x(g,g_{k-2})\}$.
\end{enumerate}

\medskip

\noindent To prove Theorem \ref{igualdad polinomios} let $k\in \{1,\ldots,h+1\}$ be the smallest integer such that $d_x(f,g)\leq d_x(f,f_{k-1})$. Then for $k>1$ we have $d_x(f,g)> d_x(f,f_{k-2})$ and by Proposition \ref{ppp} $k\leq h'+1$ and (a') holds. 
\medskip

\noindent By Proposition \ref{rrr} we get $d_x(f,g)\leq d_x(g,g_{k-1})$.
\medskip

\noindent To check (c') suppose that $d_x(f,g)<d_x(f,f_{k-1})$ and $d_x(f,g)<d_x(g,g_{k-1})$. Then by the STI $d_x(g,f_{k-1})=d_x(f,g)$ and $d_x(f,g_{k-1})=d_x(f,g)$. 
Thus $d_x(f,g)=d_x(g,f_{k-1})<d_x(g,g_{k-1})$ and using again
Áthe STI to the power series $g, f_{k-1}$ and $g_{k-1}$ we get 
\[
d_x(f_{k-1},g_{k-1})=\inf \{d_x(g,f_{k-1}), d_x(g,g_{k-1})\}=d_x(g,f_{k-1})=d_x(f,g)
\]

\noindent  which proves (c').

\medskip

\noindent Suppose that $k>1$. Then $d_x(f,g)>d_x(f,f_{k-2})$ by the definition of $k$ and $d_x(g,g_{k-2})=d_x(f,f_{k-2})$ by (a'). This proves (d'). The assertion on the key polynomials follows from Proposition \ref{ppp}.
\hfill $\blacksquare$

\section{Application to polynomial automorphisms}
\label{Kulk}

\noindent  In \cite{van der Kulk} van der Kulk proved a theorem on polynomial automorphisms of the plane generalizing a previous result of 
Jung \cite{Jung} to  the case of arbitrary characteristic. The proof of van der Kulk is based on a lemma on the intersection multiplicity of branches proved using the Hamburger-Noether expansions (see also \cite[Remark 6.3.1, p. 67]{Russell}).

\medskip

\noindent As application of our main result we prove here a property of  intersection multiplicities of branches which implies van der Kulk's lemma 
(see \cite{Ploski} for char $\mathbf{K}=0$).

\begin{prop}\label{congruencia} Let $\{f=0\}$ and $\{g=0\}$ be two different
branches and let $\{l=0\}$ be a smooth branch. Suppose that
$n=i_0(f,l)<+\infty$, $n'=i_0(g,l)<+\infty$ and let $d=\gcd(n,n')$. Then $i_0(f,g)\equiv 0$ mod
$\left(\frac{n}{d}\;\hbox{\rm or }\frac{n'}{d}\right)$.
\end{prop}

\noindent \begin{proof}
We may assume that $n,n'>1$ and $l=x$. Let $k>0$ be the integer as in  Corollary \ref{c:1}. Then we have $i_0(f,g)\leq \inf \{e'_{k-1}v_k,e_{k-1}v'_k\}$. We claim that 

\begin{equation}
\label{i2} i_0(f,g)\equiv 0 \;\;\hbox{\rm mod }(e_{k-1}\;\hbox{\rm
or }e'_{k-1}).
\end{equation}

\noindent In fact, it is clear when the equality $i_0(f,g)=\inf \{e'_{k-1}v_k,e_{k-1}v'_k\}$ holds. If $i_0(f,g)<\inf \{e'_{k-1}v_k,e_{k-1}v'_k\}$ then we conclude (\ref{i2}) from  the second part of Corollary \ref{c:1}.

\medskip

\noindent By the first part of Corollary \ref{c:1}  and Remark \ref{nota} we get 
\begin{equation}
\label{i1} \frac{n}{e_i}=\frac{n'}{e'_i} \;\;\hbox{\rm for }i<k.
\end{equation}

\noindent From (\ref{i1}) we have

\begin{equation}
\label{i3} e_{k-1}\equiv 0 \;\;\left(\hbox{\rm mod }\frac{n}{d}\right)\;\;and
\;\; e'_{k-1}\equiv 0 \;\;\left(\hbox{\rm mod }\frac{n'}{d}\right).
\end{equation}

\noindent Now (\ref{i2}) and (\ref{i3}) imply the proposition.
\end{proof}

\medskip

\noindent Using Proposition \ref{congruencia} we will prove the
following basic property of polynomial automorphisms of the plane.

\begin{teorema}( \cite [Lemma on page 36]{van der Kulk})
\label{van der Kulk}
Let the mapping $(P,Q):\bK^2\longrightarrow \bK^2$ be a polynomial
automorphism. Then of the two integers $m=\deg P$,
$n=\deg Q$ one divides the other.
\end{teorema}

\noindent \begin{proof}
Let $C$ and $D$ be projective curves with affine equations $P=0$ and
$Q=0$ respectively. Then $\deg D=n$, $\deg C=m$ and each of the
curves $C,D$ has exactly one branch at infinity (see  \cite [p. 37]{van der Kulk}).
 By B\'ezout's Theorem
these branches intersect with multiplicity $i=mn-1$. The line at
infinity cuts the branches of $C$ and $D$ with multiplicities $m$
and $n$ respectively. Thus by Proposition \ref{congruencia} we get
$i\equiv 0\;\;(\hbox{\rm mod }\frac{m}{d}\;\hbox{\rm or }
\frac{n}{d})$,  where $d=\gcd(m,n)$. This implies that
$m$ divides $n$ or $n$ divides $m$, since $i=mn-1$.
\end{proof}

\section{Intersection numbers of two branches}
\label{sBayer}
\noindent Let $(v_0, \ldots, v_h)$ and $(v'_0, \ldots, v'_{h'})$ be two characteristic sequences.
Put $e_i=\gcd(v_0,\ldots,v_i)$, $e'_j=\gcd(v'_0,\ldots,v'_j)$ for $0\leq i\leq h$ and $0\leq j\leq h'$. By convention  $v_{h+1}=v'_{h'+1}=+\infty$ and $e_{-1}=e'_{-1}=0$.
 Let 
 \[
 \rho:=\max \left\{ i\in \bN\;:\; \frac{v_j}{v_0}=\frac{v'_j}{v'_0},\;\hbox{\rm for } j\leq i, \; i\leq \min(h,h')  \right \}
 \] 
 
 \noindent and $I_k:=\inf \{e_{k-1}v'_k, e'_{k-1}v_k\}$ for $k=1, \ldots, \rho+1$.  Let $I_0=0$. Observe that $e_{k-1}v'_k=e'_{k-1}v_k$ for 
$0\leq k\leq \rho$  and $I_0<I_1<\cdots <I_{\rho+1}$. We put $ \bN^+=\{N\in \bN\;:\;N>0\}$. 

\medskip

\noindent The following theorem is a strong version of Bayer's result \cite[Theorem 5]{Bayer} proved by Hefez \cite{Hefez} (see the proof of Theorem 8.5, pp. 116-117) in characteristic zero.

\begin{teorema}
\label{Bayer}
Let $\{f=0\}$ be a branch such that $\overline{\hbox{\rm char}}_x f=(v_0, \ldots, v_h)$. Let ${\mathcal B}$ be the set of branches $\{g=0\}$ such that 
$\overline{\hbox{\rm char}}_x g=(v'_0, \ldots, v'_{h'}) $ and $i_0(f,g)\neq +\infty$. Then 
\[
\{i_0(f,g)\,:\,\{g=0\}\in {\mathcal B} \}=\bigcup_{k=1}^{\rho+1} \{N\in \bN^+\,:\,I_{k-1}\leq N<I_k \;\hbox{\rm and } N\equiv 0 \, (\hbox{\rm mod }e_{k-1}e'_{k-1})\}.
\]
\end{teorema}

\begin{coro}
\label{Umea}
Let $\{f=0\}$ be a branch such that $\overline{\hbox{\rm char}}_x f=(v_0, \ldots, v_h)$.  Then

\noindent $\{i_0(f,g)\,:\, \overline{\hbox{\rm char}}_x f=\overline{\hbox{\rm char}}_x g\;\hbox{\rm and }i_0(f,g)\neq +\infty \}$\\
$=\bigcup_{k=1}^{h+1} \{N\in \bN^+\,:\,e_{k-2}v_{k-1}\leq N< e_{k-1}v_{k}  \;\hbox{\rm and } N\equiv 0 \, (\hbox{\rm mod }e_{k-1}^2)\}$.
\end{coro}

\begin{coro}
Let $\{f=0\}$ be a branch such that $\overline{\hbox{\rm char}}_x f=(v_0, \ldots, v_h)$. Then $N_0=e_{h-1}v_h$ is the smallest natural number such that 
for any $N\in \bN$ and $N\geq N_0$  there exists an irreducible power series $g\in \bK[[x,y]]$ such that $i_0(f,g)=N$.
\end{coro}

\noindent To prove Theorem \ref{Bayer} we need two lemmas.

\begin{lema}
\label{l:1}
For any integer $k$, with $1< k \leq \rho+1$, there exists an irreducible power series $g\in \bK[[x,y]]$ such that $i_0(f,g)=I_{k-1}$.
\end{lema}
\noindent \begin{proof}

\noindent Suppose that $k>1$ and let $f_0,\ldots, f_{k-2}$ be key polynomials of $f$. Let $g_i=f_i$ for $i=0,\ldots,k-2$ and $g_i=g_{i-1}^{n'_i}+\xi_i x^{a_{i,0}}g_0^{a_{i,1}}\cdots g_{i-2}^{a_{i,i-1}}$, for $i=k-1,\ldots,h'$, where $\xi_i\in \bK\backslash \{0\}$, $n'_i=\frac{e'_{i-1}}{e'_i}$ and $a_{i,0}v'_0+\cdots + a_{i,i-1}v'_{i-1}=n'_{i}v'_{i}$. Then $g_0,\ldots, g_{h'}$ are key polynomials of 
$g:=g_{h'}$ and $\overline{\hbox{\rm char}}_x f_i=\left(\frac{v_0}{e_i}, \ldots,\frac{v_i}{e_i}\right)=\left(\frac{v'_0}{e'_i}, \ldots,\frac{v'_i}{e'_i}\right)$ for $i=0,\ldots,k-2$. 

\medskip

\noindent We have

\begin{eqnarray*}
i_0(f,x^{a_{k-1,0}}g_0^{a_{k-1,1}}\cdots g_{k-3}^{a_{k-1,k-2}})&=&a_{k-1,0}i_0(f,x)+a_{k-1,1}i_0(f,g_0)+\cdots + a_{k-1,k-2}i_0(f,g_{k-3})\\
&=&a_{k-1,0}i_0(f,x)+a_{k-1,1}i_0(f,f_0)+\cdots + a_{k-1,k-2}i_0(f,f_{k-3})\\
&=& a_{k-1,0}v_0+a_{k-1,1}v_1+\cdots + a_{k-1,k-2}v_{k-2}\\
&=& e_{k-2}\left( a_{k-1,0}\frac{v_0}{e_{k-2}}+a_{k-1,1}\frac{v_1}{e_{k-2}}+\cdots + a_{k-1,k-2}\frac{v_{k-2}}{e_{k-2}}\right)\\
&= & e_{k-2}\left( a_{k-1,0}\frac{v'_0}{e'_{k-2}}+a_{k-1,1}\frac{v'_1}{e'_{k-2}}+\cdots + a_{k-1,k-2}\frac{v'_{k-2}}{e'_{k-2}}\right)\\
&=& \frac{e_{k-2}}{e'_{k-2}}n'_{k-1}v'_{k-1}=\frac{n'_{k-1}}{e'_{k-2}}e_{k-2}v'_{k-1}=\frac{n'_{k-1}}{e'_{k-2}}e'_{k-2}v_{k-1}\\
&=&n'_{k-1}v_{k-1}=n'_{k-1} i_0(f,f_{k-2})=i_0(f,f_{k-2}^{n'_{k-1}})=i_0(f,g_{k-2}^{n'_{k-1}}).
\end{eqnarray*}

\noindent Suppose that  $\xi_{k-2}$ is generic. Then $i_0(f,g_{k-1})=i_0(f,g_{k-2}^{n'_{k-1}})=n'_{k-1}v_{k-1}$. 
We get
\[
d_x(f,g_{k-1})=\frac{n'_{k-1}v_{k-1}}{v_0\left(\frac{v'_0}{e'_{k-1}}\right)}=\frac{e'_{k-1}n'_{k-1}v_{k-1}}{v_0v'_0}=\frac{e'_{k-2}v_{k-1}}{v_0v'_0}=\frac{I_{k-1}}{v_0v'_0},
\]

\noindent and 

\[
d_x(g,g_{k-1})=\frac{v'_k}{v'_0\left(\frac{v'_0}{e'_{k-1}}\right)}=\frac{e'_{k-1}v'_{k}}{v'_0v'_0}=\frac{e_{k-1}v'_{k}}{v_0v'_0}=\frac{I_{k}}{v_0v'_0}.
\]

\noindent Therefore $d_x(f,g_{k-1})<d_x(g,g_{k-1})$, thus  $d_x(f,g)=d_x(f,g_{k-1})$, that is $\frac{i_0(f,g)}{v_0v'_0}=\frac{I_{k-1}}{v_0v'_0}$ and $i_0(f,g)=I_{k-1}$.
\end{proof}

\begin{lema}
\label{l:2}
Let $\{f=0\}$ be a branch with $\overline{\hbox{\rm char}}_x f=(v_0, \ldots, v_h)$. Let $N>0$ be an integer number such that $N>e_{h-1}v_h$. Then 
there exists a branch $\{g=0\}$ such that $i_0(f,g)=N$ and  $\overline{\hbox{\rm char}}_x g= \overline{\hbox{\rm char}}_x f$. Moreover $f_0,\ldots, f_{h-1}$ are key polynomials of $g$.
\end{lema}
\noindent \begin{proof}
Let $c$ be the conductor of $\Gamma(f)$. Then $c=(n_1-1)v_1+\cdots +(n_h-1)v_h-v_0+1=(n_1v_1-v_1)+\cdots +(n_hv_h-v_h)-v_0+1<(v_2-v_1)+(v_3-v_2)+\cdots+(e_{h-1}v_h-v_h)-v_0+1=e_{h-1}v_h-v_0-v_1+1<e_{h-1}v_h$ and we may write $N=a_0v_0+\cdots+a_hv_h$, where $a_i\in \bN$ with $a_0>0$. Let $f_k$ be a key polynomial of $f$ for $k=0,1,\ldots,h-1$. Put  $g:=f+x^{a_0}f_0^{a_1}\cdots f_{h-1}^{a_{h}}$. Then $i_0(x,g)=i_0(x,f)$ since $a_0>0$, and 
$i_0(f,g)=i_0(f,x^{a_0}f_0^{a_1}\cdots f_{h-1}^{a_{h}})=a_0v_0+ \cdots+a_hv_h=N$. By the Abhyankar-Moh irreducibility criterion (see Section \ref{AM-irr}), $g$ is irreducible and  $\overline{\hbox{\rm char}}_x g= \overline{\hbox{\rm char}}_x f$.\\ 

\noindent Observe that $d_x(f,g)=\frac{N}{v_0^2}>\frac{e_{h-1}v_h}{v_0^2}\geq \frac{e_{k}v_{k+1}}{v_0^2}=d_x(f,f_k)$. Thus, by the STI, we have $d_x(f_k,g)=d_x(f_k,f),$ which implies $i_0(f_k,g)=i_0(f_k,f)=v_{k+1}$. Therefore $f_k$ is a key polynomial of $g$ for $k=0,\ldots, h-1$.
\end{proof}

\medskip

\noindent {\bf Proof of Theorem \ref{Bayer}} The inclusion $"\subset"$ follows from Corollary \ref{c:1}. Let $N>0$ be an integer such that $I_{k-1}\leq N<I_k$ and $N\equiv 0$ (mod $e_{k-1}e'_{k-1}$), where $1\leq k\leq \rho+1$. We have to prove that there exists an irreducible power series $g\in \bK[[x,y]]$ such that $\overline{\hbox{\rm char}}_x g=(v'_0, \ldots, v'_{h'}) $ and $i_0(f,g)=N$. If $N=I_{k-1}$, then the theorem follows from Lemma \ref{l:1}. Suppose that $I_{k-1}<N<I_k$. Since $N\equiv 0$ (mod $e_{k-1}e'_{k-1}$) we may write $N=e_{k-1}e'_{k-1}N_{k-1}$ for some $N_{k-1}\in \bN$.  Let $f_0,\ldots, f_{k-1}$ be key polynomials of $f$, where $\overline{\hbox{\rm char}}_x f_i=\left(\frac{v_0}{e_i}, \ldots,\frac{v_i}{e_i}\right)=\left(\frac{v'_0}{e'_i}, \ldots,\frac{v'_i}{e'_i}\right)$ for $i=0,\ldots,k-1$ by the definition of $\rho$. Put $g_i:=f_i$ for $i=0,\ldots,k-2$.\\

\noindent {\bf Claim 1:} There exists an  irreducible power series $g_{k-1}\in \bK[[x,y]]$ such that $\overline{\hbox{\rm char}}_x g_{k-1}=\overline{\hbox{\rm char}}_x f_{k-1}$ and $i_0(f_{k-1},g_{k-1})=N_{k-1}$.

\noindent Claim 1 follows from
\begin{eqnarray*}
\gcd \left(\frac{v_0}{e_{k-1}}, \ldots,\frac{v_{k-2}}{e_{k-1}}\right)\frac{v_{k-1}}{e_{k-1}}&=&\frac{e_{k-2}v_{k-1}}{e_{k-1}^2}=\frac{e'_{k-2}v_{k-1}}{e'_{k-1}e_{k-1}}=\frac{I_{k-1}}{e_{k-1}e'_{k-1}}\\
&<&\frac{N}{e_{k-1}e'_{k-1}}=N_{k-1},
\end{eqnarray*}

\noindent and Lemma \ref{l:2} applied to $\{f_{k-1}=0\}$.\\

\noindent {\bf Claim 2:} $i_0(f,g_{k-1})=e_{k-1}N_{k-1}$.\\

\noindent If $v_{\rho+1}=+\infty$ then $f_{\rho}$ is the distinguished polynomial associated with $f$. Therefore $i_0(f,g_{\rho})=i_0(f_{\rho},g_{\rho})=N_{\rho}$ ($=e_{\rho}N_{\rho}$ since $e_{\rho}=1$) by Claim 1.

\noindent Assume that $v_{\rho+1}\neq +\infty$. Then $d_x(f_{k-1},g_{k-1})=\frac{N_{k-1}}
{\left(\frac{v_0}{e_{k-1}}\right)\left(\frac{v'_0}{e'_{k-1}}\right)}=\frac{N}{v_0v'_0}<\frac{I_k}{v_0v'_0}\leq \frac{e'_{k-1}v_k}{v_0v'_0}=\frac{e_{k-1}v_k}{v_0^2}=
d_x(f,f_{k-1})$. Therefore by the STI we have $d_x(f,g_{k-1})=d_x(f_{k-1},g_{k-1})=\frac{N}{v_0v'_0}$ and $\frac{i_0(f,g_{k-1})}{v_0\left(\frac{v'_0}{e'_{k-1}}\right)}=\frac{e_{k-1}e'_{k-1}N_{k-1}}{v_0v'_0}$, which implies $i_0(f,g_{k-1})=e_{k-1}N_{k-1}$.\\

\noindent Let us finish now the proof of Theorem \ref{Bayer}. Let $g=g_{h'}$. We will check that $d_x(f,g)=d_x(f,g_{k-1})$. Firstly suppose that $v'_{\rho+1}=+\infty$. Then $g_{\rho}$ is the distinguished polynomial associated with $g$ and $d_x(f,g)=d_x(f,g_{\rho})$. If  $v'_{\rho+1}\neq +\infty$ then we have $d_x(f,g_{k-1})=\frac{e_{k-1}N_{k-1}}{v_0\left(\frac{v'_0}{e'_{k-1}}\right)}=
\frac{N}{v_0v'_0}$, $d_x(g,g_{k-1})=\frac{v'_k}{v'_0\left(\frac{v'_0}{e'_{k-1}}\right)}=\frac{e'_{k-1}v'_k}{v'_0v'_0}=\frac{e_{k-1}v'_k}{v_0v'_0}\geq\frac{I_k}{v_0v'_0}>\frac{N}{v_0v'_0}=d_x(f,g_{k-1})$ and by the STI $d_x(f,g)=d_x(f,g_{k-1})$. Therefore $\frac{i_0(f,g)}{v_0v'_0}=\frac{N}{v_0v'_0}$ and we get $i_0(f,g)=N$.
 {\small $\blacksquare$}
 
 \section{A property of the logarithmic distance}
 \label{log dist}
 \noindent Recall that the logarithmic distance $d(f,g)$ among two branches $\{f=0\}$ and $\{g=0\}$ is given by

\[
d(f,g)=\frac{i_0(f,g)}{\ord f \, \ord g}.
\]

\noindent Observe that $d(f,g)=d_x(f,g)$ when $\{x=0\}$ is transverse to $\{f=0\}$ and $\{g=0\}$.\\

\noindent If $\{f=0\}$ and $\{g=0\}$  have no common tangent then $d(f,g)=1$. The next theorem generalizes Theorem 2.7 in \cite{Alberich} to arbitrary characteristic. 

\begin{teorema}
\label{Carr}
Let $f\in \bK[[x,y]]$ be an irreducible power series and let $R>1$ be a rational number. Then there exists an irreducible power series $g\in \bK[[x,y]]$ such that $d(f,g)=R$.
\end{teorema}

\noindent \begin{proof}
Let $\overline{\hbox{\rm char}}_x f=(v_0, \ldots, v_h)$, where $v_0<v_1$. Fix a rational number $R>1$. We distinguish two cases:\\

\noindent {\bf Case 1:} There exists an integer $k$, $1\leq k\leq h$ such that $R=\frac{e_{k-1}v_k}{v_0^2}$. Then for a $(k-1)$-th key polynomial  $f_{k-1}$ of $f$ we have $i_0(f,f_{k-1})=v_k$, $\ord f=v_0$, $\ord f_{k-1}=\frac{v_0}{e_{k-1}}$ and $d(f,f_{k-1})=R$.\\

\noindent {\bf Case 2:} The number $R\neq \frac{e_{l-1}v_l}{v_0^2}$ for $l=1,\ldots, h$. Then there exists a unique $k$, $1\leq k\leq h+1$ such that
$\frac{e_{k-2}v_{k-1}}{v_0^2}< R < \frac{e_{k-1}v_k}{v_0^2}$ (recall that $e_{-1}=0$). Write $R=\frac{r}{(\frac{v_0}{e_{k-1}})^2s}$, where $\gcd(r,s)=1$. Let $s>1$. Put $(v'_0,\ldots,v'_k)=
\left( s\frac{v_0}{e_{k-1}}, \ldots, s\frac{v_{k-1}}{e_{k-1}},r\right)$. We check that

\begin{enumerate}
\item $(v'_0,\ldots,v'_k)$ is a characteristic sequence, 
\item $\frac{v'_1}{e'_0}=\frac{v_1}{e_0},\ldots, \frac{v'_{k-1}}{e'_{0}}=\frac{v_{k-1}}{e_{0}}$ and $\frac{v'_{k}}{e'_{0}}<\frac{v_{k}}{e_{0}}$, where $e'_i:=\gcd(v'_0,\ldots,v'_i)$.
\end{enumerate}

\noindent By Theorem \ref{Bayer} there exists an irreducible $g$ such that $\overline{\hbox{\rm char}}_x g=(v'_0, \ldots, v'_k)$ and $i_0(f,g)=\inf \{e'_{k-1}v_k, e_{k-1}v'_k\}$. Therefore we get $d(f,g)=\frac{i_0(f,g)}{v_0v'_0}=\inf \left\{ \frac{e'_{k-1}v_{k}}{e'_0e_{0}}, \frac{e_{k-1}v'_{k}}{e_0e'_{0}}  \right\}=\frac{e_{k-1}v'_{k}}{e_0e'_{0}}=R$. 

\medskip

\noindent Now let $s=1$. Then $R=\frac{re_{k-1}^2}{v_0^2}$ and $e_{k-2}v_{k-1}<re_{k-1}^2<e_{k-1}v_k$. By Corollary \ref{Umea} there exists an irreducible power series $g$ such that $\ord g=\ord f=v_0$ and $i_0(f,g)=re_{k-1}^2$. Clearly $d(f,g)=R$.
\end{proof}

\noindent {\small  Evelia Rosa Garc\'{\i}a Barroso\\
Departamento de Matem\'aticas, Estad\'{\i}stica e I.O. \\
Secci\'on de Matem\'aticas, Universidad de La Laguna\\
Apartado de Correos 456\\
38200 La Laguna, Tenerife, Espa\~na\\
e-mail: ergarcia@ull.es}

\medskip

\noindent {\small Arkadiusz P\l oski\\
Department of Mathematics and Physics\\
Kielce University of Technology\\
Al. 1000 L PP7\\
25-314 Kielce, Poland\\
e-mail: matap@tu.kielce.pl}

\end{document}